\documentclass[10pt]{article}

\usepackage{amsmath}
\usepackage{amsfonts}
\usepackage{amssymb}
\usepackage{framed}
\usepackage{mathtools}
\usepackage{enumerate}
\usepackage{mathrsfs}
\usepackage[hidelinks]{hyperref}
\usepackage{enumitem}
\usepackage{mathrsfs}
\usepackage{scrextend}
\usepackage{amsthm}
\usepackage{bbm}
\usepackage{thmtools}
\usepackage{thm-restate}
\usepackage{mathabx}
\usepackage[margin=1.5in]{geometry}
\usepackage{rotating}

\usepackage{xcolor}
\hypersetup{
    colorlinks,
    linkcolor={red!50!black},
    citecolor={blue!50!black},
    urlcolor={blue!80!black}
}

\usepackage{tikz-cd}

\theoremstyle{plain}
\newtheorem{thm}{Theorem}[section]
\newtheorem{lem}[thm]{Lemma}
\newtheorem{prop}[thm]{Proposition}

\theoremstyle{definition}
\newtheorem{defn}[thm]{Definition}

\theoremstyle{remark}
\newtheorem*{rem}{Remark}

\newtheorem*{notn}{Notation}

\tikzset{
  symbol/.style={
    draw=none,
    every to/.append style={
      edge node={node [sloped, allow upside down, auto=false]{$#1$}}}
  }
}

\newcommand{\Spec}{\textrm{Spec} \hspace{0.15em} }

\newcommand\restr[2]{{
	\left.\kern-\nulldelimiterspace
	#1
	\vphantom{\big|}
	\right|_{#2}
	}}
\newcommand{\an}[1]{#1^{\textrm{an}}}

\newcommand{\Set}{\textrm{Set}}
\newcommand{\Sch}{\textrm{Sch}}

\newcommand{\Hom}{\textrm{Hom}}

\newcommand{\ch}[1]{\widecheck{{#1}}}

\newcommand{\codim}{\textrm{codim}}

\newcommand{\GL}{\textrm{GL}}

\usepackage{amsmath,calligra,mathrsfs}

\usepackage[cal=boondoxo]{mathalfa}

\DeclareMathOperator{\sheafhom}{\mathcal{H \kern -1pt o \kern -2pt m}}
\DeclareMathOperator{\sheafend}{\mathcal{E \kern -1pt n \kern -2pt d}}

\title{Effective Methods for Diophantine Finiteness}
\author{David Urbanik}

\begin{document}

\maketitle

\begin{abstract}
Let $K \subset \mathbb{C}$ be a number field, and let $\mathcal{O}_{K,N} = \mathcal{O}_{K}[N^{-1}]$ be its ring of $N$-integers. Recently, Lawrence and Venkatesh proposed a general strategy for proving the Shafarevich conjecture for the fibres of a smooth projective family $f : X \to S$ defined over $\mathcal{O}_{K,N}$. To carry out their strategy, one needs to be able to decide whether the algebraic monodromy group $\mathbf{H}_{Z}$ of any positive-dimensional geometrically irreducible subvariety $Z \subset S_{\mathbb{C}}$ is ``large enough'', in the sense that a certain orbit of $\mathbf{H}_{Z}$ in a variety of Hodge flags has dimension bounded from below by a certain quantity. In this article we give an effective method for deciding this question. Combined with the effective methods of Lawrence-Venkatesh for understanding semisimplifications of global Galois representations using $p$-adic Hodge theory, this gives a fully effective strategy for solving Shafarevich-type problems for arbitrary families $f$. 
\end{abstract}

\tableofcontents

\section{Introduction}

Let $f : X \to S$ be a smooth projective family defined over $\mathcal{O}_{K, N} = \mathcal{O}_{K}[N^{-1}]$, with $K \subset \mathbb{C}$ a number field, $\mathcal{O}_{K}$ its ring of integers, and $N$ any non-zero element of $\mathbb{Z}$. In a recent paper \cite{LV}, Lawrence and Venkatesh proposed a strategy for bounding the dimension of the Zariski closure $\overline{S(\mathcal{O}_{K,N})}^{\textrm{Zar}}$ of $S(\mathcal{O}_{K,N})$. Since for any point $s \in S(\mathcal{O}_{K,N})$ the fibre $X_{s}$ has good reduction away from $N$, the problem of bounding the dimension of $\overline{S(\mathcal{O}_{K,N})}^{\textrm{Zar}}$ can be interpreted as Shafarevich-type problem for the family $f$. In particular, the Shaferevich conjectures for curves, abelian varieties, and K3 surfaces are equivalent to showing that $\dim \overline{S(\mathcal{O}_{K,N})}^{\textrm{Zar}} = 0$ for particular choices of $f$.

The Lawrence-Venkatesh strategy has been implemented to reprove various classical diophantine finiteness results, including the Mordell conjecture \cite{LV}, the finiteness of $S$-units \cite{LV}, and Seigel's theorem \cite{noordman2021siegels}. However it is of particular interest to apply the strategy in situations where $\dim S > 1$. In this situation the problem becomes much harder, because one needs some understanding of the monodromy of the family $f$ over essentially arbitrary geometrically irreducible subvarieties $Z \subset S_{\mathbb{C}}$. To explain what we mean, let $\mathbb{V} = R^{i} f_{*} \mathbb{Z}$ and define for each such $Z \subset S_{\mathbb{C}}$:
\begin{defn}
The algebraic monodromy group $\mathbf{H}_{Z}$ of $Z$ is the identity component of Zariski closure of the monodromy representation associated to $\restr{\mathbb{V}}{Z^{\textrm{nor}}}$, where $Z^{\textrm{nor}} \to Z$ is the normalization.
\end{defn}
\noindent The variation $\mathbb{V}$ determines a flag variety $\ch{L}$, on which the group $\mathbf{H}_{Z}$ can be said to act after identifying $\ch{L}$ with the variety of Hodge flags on some fibre $\mathbb{V}_{s}$ for $s \in Z(\mathbb{C})$. Let $\varphi : S \to \Gamma \backslash D$ be a period map determined by $\mathbb{V}$ with $D \subset \ch{L}$ the complex submanifold of polarized Hodge flags. Then the key quantity relevant for the Lawrence-Venkatesh method is
\begin{equation}
\Delta = \min_{Z, s} \left[ \dim (\mathbf{H}_{Z} \cdot F^{\bullet}_{s}) - \dim \varphi(Z) \right] , 
\end{equation}
where the minimum is taken over all positive-dimensional geometrically irreducible subvarieties $Z \subset S_{\mathbb{C}}$, points $s \in Z(\mathbb{C})$, and where $F^{\bullet}_{s}$ is the Hodge flag on $\mathbb{V}_{s}$.\footnote{We not here that $\varphi(Z)$ is analytically constructible, for instance by applying the main result of \cite{OMINGAGA}, so its dimension makes sense.} In particular, the Lawrence-Venkatesh method produces an integer $k$, and shows that if $\Delta \geq k$, then the Shafarevich conjecture holds for $f$.

For successful applications of the Lawrence-Venkatesh strategy for the Shafarevich problem in situations when $\dim S > 1$ we know of only the paper \cite{lawrence2020shafarevich} of Lawrence and Sawin, who are able to apply this strategy beyond the first induction step to prove a Shafarevich conjecture for hypersurfaces lying inside a fixed abelian variety $A$. Their methods require the auxilliary use of a Tannakian category associated to $A$, and it seems unclear what to do without this abelian variety structure present.

Our main result is as follows:

\begin{thm}
\label{mainthm}
Consider a smooth projective family $f : X \to S$ defined over $\mathcal{O}_{K,N}$ and with $S$ smooth and quasi-projective, and given an integer $d$, define
\begin{equation}
\Delta_{d} = \min_{Z, s} \left[ \dim (\mathbf{H}_{Z} \cdot F^{\bullet}_{s}) - \dim \varphi(Z) \right] , 
\end{equation}
where the minimum ranges over all geometrically-irreducible subvarieties $Z \subset S_{\mathbb{C}}$ with $\dim Z > d$ and points $s \in Z(\mathbb{C})$. Then there exists an effective procedure which outputs an infinite sequence of integers
\[ \kappa(1) \leq \kappa(2) \leq \cdots \leq \kappa(r) \leq \cdots < \Delta_{d} \]
such that for some $r = r_{0}$ we have $\kappa(r_{0}) = \Delta_{d} - 1$.

If moreover the period map $\varphi$ is quasi-finite, one can determine $r_{0}$.
\end{thm}

\begin{rem}
Let us make absolutely clear what is meant by the term ``effective procedure'' in \autoref{mainthm}. We mean that there exists an infinite-time algorithm (for instance, a non-halting Turing machine), which outputs a sequence of integers $\{ \kappa(r) \}_{r = 1}^{\infty}$, with the integer $\kappa(r)$ being outputted at some finite time $t_{r}$ depending on $r$. Moreover, one also has at time $t_{r}$ a proof that $\kappa(r) < \Delta_{d}$. Therefore, after time $t_{r}$, one can stop the algorithm and use the bound $\kappa(r) < \Delta_{d}$ as input to the Lawrence-Venkatesh method. One is also guaranteed that there is some $r_{0}$ so that at time $t_{r_{0}}$ the bound $\kappa(r_{0}) < \Delta_{d}$ is best possible, however one does not necessarily have a method to determine $r_{0}$ unless $\varphi$ is quasi-finite. Finally, the algorithm can be described entirely in terms of algebro-geometric computations involving algebraically constructible sets, and implicit in the proof is a description of how to implement it.
\end{rem}

There is no fundamental obstruction which requires us to restrict to quasi-finite $\varphi$ for the second part of \autoref{mainthm}. Rather, the second part of \autoref{mainthm} references some delicate arguments in \cite{urbanik2021sets} which are only enough to handle the quasi-finite case directly, and recalling enough of the machinery of \cite{urbanik2021sets} to carry out the argument for the general case would lead us too far astray from the main ideas. We note that one does not actually need to determine the integer $r_{0}$ referenced in \autoref{mainthm} to apply the machinery of Lawrence and Venkatesh: one wants to be able to compute the best possible lower bound for $\Delta_{d}$, but one is not required to actually prove that the bound one has is optimal in order to deduce diophantine finiteness results.

\subsection{The Approach of Lawrence and Venkatesh}

We begin with a preliminary observation. To show that $\dim \overline{S(\mathcal{O}_{K,N})}^{\textrm{Zar}} \leq d$, it suffices to show, for any irreducible subscheme $T \subset S$ of dimension $> d$ and defined over $\mathcal{O}_{K,N}$, that the Zariski closure $\overline{T(\mathcal{O}_{K,N})}^{\textrm{Zar}}$ is a proper algebraic subscheme of $T$. We therefore fix such a subscheme $T \subset S$ with $\dim T > d$, and seek to show that $\dim \overline{T(\mathcal{O}_{K,N})}^{\textrm{Zar}} < \dim T$.

Fix a prime $p \in \mathbb{Z}$ not dividing $N$, and let $t \in T(\mathcal{O}_{K,N})$ be a point. It is conjectured that for any $i$ the representation of $\textrm{Gal}(\overline{K}/K)$ on $H^{i}_{\textrm{\'et}}(X_{t,\overline{K}}, \mathbb{Q}_{p})$ is semisimple. If this result were to be known for all such $t$, an argument of Faltings shows that for $t \in T(\mathcal{O}_{K,N})$ there are at most finitely many possibilities for the isomorphism class of the representation of $\textrm{Gal}(\overline{K}/K)$ on $H^{i}_{\textrm{\'et}}(X_{t,\overline{K}}, \mathbb{Q}_{p})$. To establish the non Zariski-density of $T(\mathcal{O}_{K,N})$ it would then suffice to show that the fibres of the map
\[ t \in T(\mathcal{O}_{K,N})  \hspace{1em} \xrightarrow{\tau} \hspace{1em} \big\{ \textrm{Gal}(\overline{K}/K)\textrm{-rep. on }H^{i}_{\textrm{\'et}}(X_{t,\overline{K}}, \mathbb{Q}_{p}) \big\}\hspace{0.5em} \big/ \hspace{0.5em} \textrm{iso}.  \]
are not Zariski dense. As explained by Lawrence and Venkatesh in \cite{LV}, this is essentially the original argument for the Mordell conjecture due to Faltings.

The problem with applying this strategy in general is twofold. First, the semisimplicity of $H^{i}_{\textrm{\'et}}(X_{t,\overline{K}}, \mathbb{Q}_{p})$ is not known, and for most choices of $f$ appears out of reach. Secondly, without a good geometric interpretation of the \'etale cohomology $H^{i}_{\textrm{\'et}}(X_{t,\overline{K}}, \mathbb{Q}_{p})$ it is difficult to understand $\tau$. The key insight in the paper of Lawrence and Venkatesh is that one may potentially overcome both problems by passing to a $p$-adic setting where they are more managable.

Instead of considering the global Galois representation $\rho_{t} : \textrm{Gal}(\overline{K}/K) \to H^{i}_{\textrm{\'et}}(X_{t,\overline{K}}, \mathbb{Q}_{p})$, we consider its semisimplification $\rho^{\textrm{ss}}_{t}$, and restrict $\rho_{t}$ along the map $\textrm{Gal}(\overline{K_{v}}/K_{v}) \to \textrm{Gal}(\overline{K}/K)$ induced by a fixed embedding $\overline{K} \hookrightarrow \overline{K_{v}}$ to obtain $\rho_{t,v}$, where $v$ is a fixed place above $p$. The functors of $p$-adic Hodge theory tell us that the representation $\rho_{t,v}$ determines a triple $(H^{i}_{\textrm{dR}}(X_{t,K_{v}}), \phi_{t}, F^{\bullet}_{t})$, where $F^{\bullet}_{t}$ is the Hodge filtration on de Rham cohomology and $\phi_{t}$ is the crystalline Frobenius at $t$. If we somehow manage to consider the data $(H^{i}_{\textrm{dR}}(X_{t,K_{v}}), \phi_{t}, F^{\bullet}_{t})$ up to ``semisimplification'' (in the sense that we identify such triples when the associated \emph{global} Galois representations have isomorphic semisimplifications), our problem is then to study the fibres of the map
\[ t \in T(\mathcal{O}_{K,N}) \hspace{0.5em} \xrightarrow{\tau_{p}} \hspace{0.5em} \big\{\textrm{``semisimplifications'' of } (H^{i}_{\textrm{dR}}(X_{t,K_{v}}), \phi_{t}, F^{\bullet}_{t}) \big\} \hspace{0.5em} \big/ \hspace{0.5em} \textrm{iso}. \]
and show they lie in a Zariski closed subscheme of smaller dimension.

Next, we make the elementary observation that to bound the dimension of the Zariski closure of $T(\mathcal{O}_{K,N})$, it suffices to cover $T(\mathcal{O}_{K,v})$ by finitely many $v$-adic disks $D_{1}, \hdots, D_{k}$ and bound the dimension of the Zariski closure of $D_{i} \cap T(\mathcal{O}_{K,N})$ for each $i$; here $\mathcal{O}_{K,v}$ is the ring of $v$-adic integers. It can then be shown that there exists such a cover for which the Hodge bundle $\mathcal{H} = R^{i} f_{*} \Omega^{\bullet}_{X/S}$ can be trivialized rigid-analytically over each $D_{i}$, moreover with respect to each such trivialilzation the Frobenius operator $\phi_{t}$ is independent of $t \in D_{i}$. The problem then reduces to studying a varying filtration $F^{\bullet}_{t}$ on a fixed vector space $V_{p} = H^{i}_{\textrm{dR}}(X_{t_{0}})$ for some $t_{0} \in D_{i}$. In particular, we obtain a rigid-analytic map
\[ D_{i} \hspace{0.5em} \xrightarrow{\psi_{p}} \hspace{0.5em} \underbrace{\{ \hspace{0.3em} \textrm{Hodge filtrations on }V_{p} \hspace{0.3em} \}}_{\ch{L}_{p}} , \]
and those points of $\ch{L}_{p}$ arising from points $t \in T(\mathcal{O}_{K,N}) \cap D_{i}$ lie inside finitely many subvarieties $O_{i1}, \hdots, O_{i\ell}$ of $\ch{L}_{p}$ corresponding to the finitely many possible isomorphism classes of $\rho^{\textrm{ss}}_{t}$.

We are now faced with the problem of understanding the intersections $\psi_{p}(D_{i}) \cap O_{ij}$, and showing that their inverse images under $\psi_{p}^{-1}$ lie in an algebraic subscheme of smaller dimension. One part of this problem is to understand the dimensions of the varieties $O_{ij}$, and to show that they are sufficiently small: this step is carried out successfully for the families of hypersurfaces studied both by Lawrence-Venkatesh and Lawrence-Sawin, and appears to be managable in general. The more difficult object to control is $\psi_{p}(D_{i})$, for which we need to understand the variation of the filtration $F^{\bullet}$ over $D_{i}$. This, in turn, is governed by the Gauss-Manin connection $\nabla : \mathcal{H} \otimes \Omega^{1}_{T} \to \mathcal{H}$, which exists universally over $S$ after possibly increasing $N$; we may adjust $p$ so that it does not divide $N$ if necessary. The fact that $\nabla$ exists universally over $\mathcal{O}_{K,N}$ means that the same system of differential equations satisfied by $\psi_{p}$ at $t \in T(\mathcal{O}_{K,N}) \cap D_{i}$ is also satisfied by a Hodge-theoretic period map $\psi : B \to \ch{L}$ on a sufficiently small analytic neighbourhood $B \subset T(\mathbb{C})$ containing $t$, where $\ch{L}$ is a variety of Hodge flags. In particular, one can prove that the Zariski closures of $\psi_{p}(D_{i})$ and $\psi(B)$ have the same dimension.

The final step, which is to show that the Zariski closure of $T(\mathcal{O}_{K,N}) \cap D_{i}$ in $T$ has smaller dimension, is completed as follows. The Ax-Schanuel Theorem \cite{AXSCHAN} for variations of Hodge structures due to Bakker-Tsimerman shows that if $V$ is an algebraic subvariety of $\ch{L}$ of dimension at most $\dim \overline{\psi(B)}^{\textrm{Zar}} - \dim \psi(B)$,\footnote{The dimension $\dim \psi(B)$ can once again, at least for open neighbourhoods $B$ with a sufficiently mild geometry, be made sense of the dimension of a locally constructible analytic set. Alternatively one can replace $\dim \psi(B)$ with $\dim \varphi(Z)$, where $\varphi$ is as before and $Z \subset T_{\mathbb{C}}$ is a component containing $B$.} then the inverse image under $\psi$ of the intersection $\psi(B) \cap V$ lies in an algebraic subvariety of $T_{\mathbb{C}}$ of smaller dimension. Choosing an isomorphism $\mathbb{C} \cong \overline{K_{v}}$ one can transfer this fact to the same statement for the map $\psi_{p}$ and in particular for $V = O_{ij}$. Our problem is finally reduced to giving a lower bound for the difference $\dim \overline{\psi(B)}^{\textrm{Zar}} - \dim \psi(B)$. Our main result now reads:

\begin{thm}
\label{mainthm2}
Define the quantity
\[ \Delta_{d} := \min_{Z, \psi} \left[ \dim \overline{\psi(B)}^{\textrm{Zar}} - \dim \psi(B) \right] , \]
where $Z$ ranges over all irreducible complex algebraic subvarieties $Z \subset S_{\mathbb{C}}$ of dimension greater than $d$, and where $\psi$ is any complex analytic period map determined by the variation of Hodge structures $\mathbb{V} = R^{i} f_{*} \mathbb{Z}$ and defined on a neighbourhood $B \subset Z(\mathbb{C})$. Then there exists an effective procedure which outputs an infinite sequence of lower bounds
\[ \kappa(1) \leq \kappa(2) \leq \cdots \leq \kappa(r) \leq \cdots < \Delta_{d} \]
such that for some $r = r_{0}$ we have $\kappa(r_{0}) = \Delta_{d} - 1$.

If moreover the period map $\varphi$ is quasi-finite, one can determine $r_{0}$.
\end{thm}

 We note that by \cite[Lem. 4.10(ii)]{urbanik2021sets} it is also a consequence of the Bakker-Tsimerman Theorem that when $Z$ is geometrically irreducible, we have $\overline{\psi(B)}^{\textrm{Zar}} = \mathbf{H}_{Z} \cdot \psi(t)$ for any point $t \in Z(\mathbb{C})$, which recovers the statement of \autoref{mainthm}. 

\subsection{Basic Idea of the Method}
\label{methodsketch}

We may observe that the computation of the bound described in \autoref{mainthm2} is a purely Hodge-theoretic problem, i.e., it concerns only properties of the integral variation $\mathbb{V} = R^{i} f_{*} \mathbb{Z}$ of Hodge structures on $S_{\mathbb{C}}$. Let $\mathcal{Q} : \mathbb{V} \otimes \mathbb{V} \to \mathbb{Z}$ be a polarization of $\mathbb{V}$, and let $(V, Q)$ be a fixed polarized lattice isomorphic to one (hence any) fibre of $(\mathbb{V}, \mathcal{Q})$; as it causes no harm, we will assume that $V = \mathbb{Z}^{m}$ for some $m$, and therefore sometimes write $\GL_{m}$ for $\GL(V)$. Let $D$ be the complex manifold parametrizing polarized Hodge structures on $(V, Q)$ with the same Hodge numbers as $(\mathbb{V}, \mathcal{Q})$. A point $h \in D$ we may view as a morphism $h : \mathbb{S} \to \GL(V)_{\mathbb{R}}$, where $\mathbb{S}$ is the Deligne torus, and the Mumford-Tate group $\textrm{MT}(h)$ is the $\mathbb{Q}$-Zariski closure of $h(\mathbb{S})$. 

To present our method, we introduce some terminology:

\begin{notn}
We denote by $\ch{L}$ the $\mathbb{Q}$-algebraic variety of all Hodge flags on the lattice $V$, not necessarily polarized. We note that $D$ is an open submanifold of a closed $\mathbb{Q}$-algebraic subvariety $\ch{D} \subset \ch{L}$.
\end{notn}

\begin{defn}
Given two subvarieties $W_{1}, W_{2} \subset \ch{L}$, we say that $W_{1} \sim_{\GL} W_{2}$ if there exists $g \in \GL_{m}(\mathbb{C})$ such that $g \cdot W_{1} = W_{2}$. Given a variety $W \subset \ch{L}$, we call the equivalence class $\mathcal{C}(W)$ under $\sim_{\GL}$ a \emph{type}. The dimension of a type $\mathcal{C}(W)$ is the dimension of $W$.
\end{defn}

\begin{defn}
We say that a type $\mathcal{C}$ is \emph{Hodge-theoretic} if $\mathcal{C} = \mathcal{C}(W)$, where $W = N(\mathbb{C}) \cdot h$ for $h \in D$ and $N$ a $\mathbb{Q}$-algebraic normal subgroup of $\textrm{MT}(h)$.
\end{defn}

\vspace{0.5em}

The first step in our algorithm is:

\begin{quote}
\textbf{Step One:} Compute a finite list of types $\mathcal{C}_{1}, \hdots, \mathcal{C}_{\ell}$ such that every Hodge-theoretic type appears somewhere in the list.
\end{quote}

\noindent When we say to compute a type $\mathcal{C}$, we mean to compute a representative $W \subset \ch{L}$ such that $\mathcal{C} = \mathcal{C}(W)$. That there are only finitely many Hodge-theoretic types is shown in \autoref{finmantypes} below.

The problem given in Step One is solved in \cite[Prop. 5.4]{urbanik2021sets}; we will say little about it here. It is related to the problem of classifying Mumford-Tate groups up to conjugacy by $\GL_{m}(\mathbb{C})$, for which one can use a constructive version of the proof in \cite[Thm. 4.14]{hodgelocivoisin}. It is also similar to the problem of classifying Mumford-Tate domains as studied in \cite[Chap. VII]{GGK}. We note that the methods of \cite[Chap. VII]{GGK}, when they can be carried out effectively, result in an approach for which $\mathcal{C}_{1}, \hdots, \mathcal{C}_{\ell}$ will be exactly the set of Hodge-theoretic types. 

\vspace{0.5em}

The second step is more involved, and is the crux of our method. To describe it we need to introduce some terminology. 

\begin{defn}
\label{locperdef}
A \emph{local period map} is a map $\psi : B \to \ch{L}$ obtained as a composition $\psi = q \circ A$, where:
\begin{itemize}
\item[(i)] The set $B \subset S(\mathbb{C})$ is a connected analytic neighbourhood on which $\mathbb{V}$ is constant and $F^{k} \mathcal{H}$ is trivial for each $k$, where $\mathcal{H} = \mathbb{V} \otimes \mathcal{O}_{\an{S}}$.
\item[(ii)] The map $A : B \to \GL_{m}(\mathbb{C})$ is a varying filtration-compatible period matrix over $B$. More precisely, there exists a basis $v^{1}, \hdots, v^{m}$ for $\mathcal{H}(B)$, compatible with the filtration in the sense that $F^{k} \mathcal{H}(B)$ is spanned by $v^{1}, \hdots, v^{i_{k}}$ for some $i_{k}$, and a flat frame $b^{1}, \hdots, b^{m}$ for $\mathbb{V}_{\mathbb{C}}(B)$, such that $A(s)$ is the change-of-basis matrix from $v^{1}_{s}, \hdots, v^{m}_{s}$ to $b^{1}_{s}, \hdots, b^{m}_{s}$.
\item[(iii)] The map $q : \GL_{m} \to \ch{L}$ sends a matrix $M$ to the Hodge flag $F_{M}^{\bullet}$ defined by the property that $F^{k}_{M}$ is spanned by the first $i_{k}$ columns.
\end{itemize}
\end{defn}

\noindent To summarize the preceding definition: a local period map is exactly a period map on $B$ except one does not necessarily compute periods with respect to the integral lattice $\mathbb{V}(B) \subset \mathbb{V}_{\mathbb{C}}(B)$ but is instead allowed to consider periods with respect to a more general complex flat frame. There is a natural $\GL_{m}(\mathbb{C})$-action on the set of germs of local period maps at a point $s \in S(\mathbb{C})$, where $M \in \GL_{m}(\mathbb{C})$ acts on the map $\psi = q \circ A$ to give $M \cdot \psi = q \circ (M \cdot A)$. This action corresponds exactly to a change of the flat frame $b^{1}, \hdots, b^{m}$, and all germs of local period maps at $s$ lie in a single $\GL_{m}(\mathbb{C})$-orbit.

The construction of a local period map $\psi : B \to \ch{L}$ involves picking a basis $b^{1}, \hdots, b^{m}$ of $\mathbb{V}_{\mathbb{C}}(B)$, and hence choosing an isomorphism $\mathbb{V}_{\mathbb{C}}(B) \simeq \mathbb{C}^{m}$. When working with a local period map, we will always assume that such a basis has been choosen, and hence identify subgroups of $\GL(\mathbb{V}_{\mathbb{C}}(B))$ with subgroups of $\GL_{m}(\mathbb{C})$. In particular, if $Z \subset S_{\mathbb{C}}$ is a geometrically irreducible subvariety which intersects $B$, we have an induced action of $\mathbf{H}_{Z}$ on $\ch{L}$.

Lastly, we need:

\begin{defn}
Given two types $\mathcal{C}_{1}$ and $\mathcal{C}_{2}$, we say that $\mathcal{C}_{1} \leq \mathcal{C}_{2}$ if there exists $W_{i} \subset \ch{L}$ for $1 = 1, 2$ such that $\mathcal{C}_{i} = \mathcal{C}(W_{i})$ and $W_{1} \subset W_{2}$.
\end{defn}

\begin{defn}
\label{vartypedef}
Given a local period map $\psi : B \to \ch{L}$ and a geometrically irreducible subvariety $Z \subset S_{\mathbb{C}}$ intersecting $B$ at $s$, we call $\mathcal{C}(\overline{\psi(B \cap Z)}^{\textrm{Zar}}) = \mathcal{C}(\mathbf{H}_{Z} \cdot \psi(s))$ the \emph{type} of $Z$, and denote it by $\mathcal{C}(Z)$.
\end{defn}

\noindent For well-definedness, see \autoref{welldeflem} below. From Step One, we have computed a finite list $\mathcal{C}_{1}, \hdots, \mathcal{C}_{\ell}$ of types containing all types that can arise from the variation $\mathbb{V}$. Our next task is then:

\begin{quote}
\textbf{Step Two:} For each type $\mathcal{C}_{i}$ appearing in the list, compute a differential system $\mathcal{T}(\mathcal{C}_{i})$ on $S$ characterized by the property that an algebraic subvariety $Z \subset S_{\mathbb{C}}$ is an integral subvariety for $\mathcal{T}(\mathcal{C}_{i})$ if and only if $\mathcal{C}(Z) \leq \mathcal{C}_{i}$, and determine the dimension of a maximal integral subvariety for this system.
\end{quote}

\noindent We explain precisely what we mean by ``differential system'' in \autoref{secthree}; actually our method does something more subtle than Step Two due to the fact that we can only approximate $\mathcal{T}(\mathcal{C}_{i})$ up to some finite order, but for expository purposes this is the essential point. After this, we will see the problem is reduced to analyzing which of the differential systems $\mathcal{T}(\mathcal{C}_{i})$ admit algebraic solutions of ``exceptional'' dimension, which can be carried out using tools from functional transcendence.

\subsection{Acknowledgements}

The author thanks Brian Lawrence, Akshay Venkatesh, and Will Sawin for comments on a draft of this manuscript.

\section{Algebraic Monodromy Orbits up to Conjugacy}

In this section we describe an effective method for solving ``Step One'' as posed in \autoref{methodsketch}. We will also prove some preliminary facts about types used in the introduction, and we continue with the notation established there. We will work in the context of a general polarizable integral variation of Hodge structure $\mathbb{V}$ on the complex algebraic variety $S$, not necessarily coming from a projective family as in the introduction.

\subsection{Basic Properties of Types}

\begin{lem}
\label{finmantypes}
For any geometrically irreducible subvariety $Z \subset S$ and any local period map $\psi : B \to \ch{L}$ with $Z \cap B$ non-empty, we have
\[ \overline{\psi(Z \cap B)}^{\textrm{Zar}} = \mathbf{H}_{Z} \cdot \psi(s) , \]
for any point $s \in Z(\mathbb{C})$.
\end{lem}

\begin{proof}
It suffices to show that 
\[ \overline{\psi(C)}^{\textrm{Zar}} = \mathbf{H}_{Z} \cdot \psi(s) , \]
for each analytic component $C \subset Z \cap B$ separately, with $s$ a point of $C$. By acting on $\psi$ by an element of $\GL_{m}(\mathbb{C})$, the claim can be reduced to the situation where the periods which determine $\psi$ are computed with respect to a basis for the integral lattice $\mathbb{V}(B)$, and then the claim follows from \cite[Lem. 4.10(ii)]{urbanik2021sets}.
\end{proof}

\begin{lem}
\label{welldeflem}
The equivalence class under $\sim_{\GL}$ of $\overline{\psi(B \cap Z)}^{\textrm{Zar}}$ is independent of $\psi$; i.e., the type of $Z \subset S$ is well-defined.
\end{lem}

\begin{proof}
Let $p : Z^{\textrm{sm}} \to Z$ be a smooth resolution, and consider the variation $p^{*} \mathbb{V}$. From the fact that germs of local period maps on $Z^{\textrm{sm}}$ with respect to the variation $p^{*} \mathbb{V}$ factor through germs of local period maps on $S$, we may reduce to the same problem for $Z^{\textrm{sm}}$ and the variation $p^{*} \mathbb{V}$, i.e., we may assume $Z = S$. By analytically continuing a fixed local period map $\psi$ to the universal covering $\widetilde{S} \to S$, we learn from the irreducibility of $\widetilde{S}$ that at each point $s \in S$, there exists a local period map $\psi_{s} : B_{s} \to \ch{L}$ such that $\overline{\psi_{s}(B_{s})}^{\textrm{Zar}} = \overline{\psi(B)}^{\textrm{Zar}}$. Since the Zariski closure of $\psi_{s}(B_{s})$ is determined by the germ of $\psi_{s}$ at $s$, and because all germs of local period maps at $s$ lie in a single $\GL_{m}(\mathbb{C})$-orbit, the result follows.
\end{proof}

\begin{lem}
\label{fintypesarise}
There are only finitely many Hodge-theoretic types.
\end{lem}

\begin{proof}
We observe that the problem reduces to the following: show they are finitely many $\GL_{m}(\mathbb{C})$-equivalence classes of pairs $(h, N)$, where
\begin{itemize}
\item[(i)] $h \in D$ is a polarized Hodge structure; and
\item[(ii)] $N$ is a $\mathbb{Q}$-algebraic connected normal subgroup of $\textrm{MT}(h)$;
\end{itemize}
where we regard $\GL_{m}(\mathbb{C})$ as acting on $h$ through its action on $\ch{L}$, and on $N$ by conjugation. Note that two such equivalent pairs will generate orbits in $\ch{L}$ equivalent under $\sim_{\GL}$. Since the groups $\textrm{MT}(h)$ are reductive and have finitely many connected normal algebraic factors, this reduces to the same problem for pairs of the form $(h, \textrm{MT}(h))$. We recall that $D$ is an open submanifold of $\ch{D}$, the flag variety of flags satisfying the first Hodge-Riemann bilinear relation (the isotropy condition), and that $\ch{D}$ is an algebraic subvariety of $\ch{L}$. We then use the fact that there are finitely many Mumford-Tate groups up to $\GL_{m}(\mathbb{C})$-conjugacy (see \cite[Thm. 4.14]{hodgelocivoisin}), and that for a fixed Mumford-Tate group $M$ the Hodge structures in $D$ with Mumford-Tate contained in $M$ lie inside finitely many $M(\mathbb{C})$-orbits in $\ch{D}$, see \cite[VI.B.9]{GGK}.
\end{proof}

\subsection{Computing Types up to Conjugacy}

In this section we give some references for carrying out Step One as described in the introduction.

\begin{prop}
\label{MTgroupequivalgo}
There exists an algorithm to compute subvarieties $W_{1}, \hdots, W_{\ell} \subset \ch{L}$ such that the set of Hodge-theoretic types is a (possibly proper) subset of $\{ \mathcal{C}(W_{1}), \hdots, \mathcal{C}(W_{\ell}) \}$.
\end{prop}

\begin{proof}
This is solved in \cite[Prop. 5.4]{urbanik2021sets}.
\end{proof}

Let us comment briefly on a different approach to Step One given in \cite[Chap. VII]{GGK}. In \cite[Chap. VII]{GGK}, the authors describe a method for classifying both Mumford-Tate groups and Mumford-Tate domains (orbits of points $h \in D$ under $\textrm{MT}(h)(\mathbb{R})$ and $\textrm{MT}(h)(\mathbb{C})$). Given an appropriate such classification, one can easily solve Step One by computing the decompositions of the groups $\textrm{MT}(h)$ that arise into $\mathbb{Q}$-simple factors. The method of \cite[Chap. VII]{GGK} is to first classify CM Hodge structures $h_{\textrm{CM}} \in D$, and then give a criterion for deciding when a Lie subalgebra of $\mathfrak{gl}(V)$ corresponds to a Mumford-Tate group generating a Mumford-Tate domain containing $h_{\textrm{CM}}$. They carry out this classification procedure successfully when $\dim V = 4$, and so for variations with Hodge numbers $(2, 2)$ and $(1, 1, 1, 1)$.

The method given for classifying CM Hodge structures given in \cite[Chap. VII]{GGK} is to observe that CM Hodge structures up to isogeny are determined by certain data associated to embeddings of CM fields, and hence the first step of the procedure in \cite[Chap. VII]{GGK} is to ``classify all CM fields of rank up to [$\dim V$] by [their] Galois group''. We are not aware of an effective method for carrying out this step.\footnote{The paper \cite{dodson1984structure} gives a potential approach by giving a method to classify certain abstract structures associated with Galois groups of CM fields. However one still needs to determine which such structures are actually associated to a concrete CM field.} It is also not clear to us precisely the sense in which the term ``classify'' is being used; i.e., we do not know what form the data of a ``classification of CM Hodge structures'' takes, and consequently what form the resulting classification of Mumford-Tate domains will have. For this reason, we were unable to apply the methods of \cite[Chap. VII]{GGK} to prove \autoref{MTgroupequivalgo}. 


\section{Differential Tools and a Jet Criterion}
\label{secthree}

In this section we introduce a collection of effectively computable algebro-geometric correspondences which can be used for studying systems of differential equations on $S$ induced by the variation $\mathbb{V}$, and then use it to solve the main problem. We have already carried out most of the work in two preceding papers \cite{periodimages} and \cite{urbanik2021sets}, so we will first need to collect some results. In this section we assume that $S$ is a $K$-variety for $K \subset \mathbb{C}$ a number field, and that $\mathbb{V}$ is a polarizable integral variation of Hodge structure on $S_{\mathbb{C}}$ such that the vector bundle $\mathcal{H} = \mathbb{V} \otimes_{\mathbb{Z}} \mathcal{O}_{\an{S_{\mathbb{C}}}}$, the filtration $F^{\bullet}$, and the connection $\nabla : \mathcal{H} \to \mathcal{H} \otimes \Omega^{1}_{S}$ all admit $K$-algebraic models. Moreover, we assume that we may effectively compute a description of these objects in terms of finitely-presented $K$-modules over an affine cover of $S$; for a justification of this assumption in the situation where $\mathbb{V}$ comes from a smooth projective $K$-algebraic family $f : X \to S$ see \cite[\S2]{urbanik2021sets}.

\subsection{The Constructive Period-Jet Correspondence}

Our algebro-geometric correspondences will be formulated using the language of \emph{jets}. Let $A^{d}_{r} = K[t_{1}, \hdots, t_{d}]/\langle t_{1}, \hdots, t_{d} \rangle^{r+1}$, and define $\mathbb{D}^{d}_{r} = \Spec A^{d}_{r}$ to be the $d$-dimensional disk of order $r$; we suppress the field $K$ in the notation. A \emph{jet space} associated to a space $X$ is a space which parametrizes maps $\mathbb{D}^{d}_{r} \to X$. More formally, for $X$ a finite-type $K$-scheme, we have:

\begin{defn}
\label{jetspacedef}
We define $J^{d}_{r} X$ to be the scheme representing the contravariant functor $\Sch_{K} \to \Set$ given by \vspace{-0.2em}
\[ T \mapsto \Hom_{K}(T \times_{K} \mathbb{D}^{d}_{r}, X), \hspace{1.5em} [T \to T'] \mapsto [\Hom_{K}(T' \times_{K} \mathbb{D}^{d}_{r}, X) \to \Hom_{K}(T \times_{K} \mathbb{D}^{d}_{r}, X)] , \]
where the natural map $\Hom_{K}(T' \times_{K} \mathbb{D}^{d}_{r}, X) \to \Hom_{K}(T \times_{K} \mathbb{D}^{d}_{r}, X)$ obtained by pulling back along $T \times_{K} \mathbb{D}^{d}_{r} \to T' \times_{K} \mathbb{D}^{d}_{r}$. 
\end{defn}

\vspace{0.5em}

\noindent That the functor defining $J^{d}_{r} X$ in \autoref{jetspacedef} is representable is handled by \cite[\S2]{periodimages}. Moreover, $J^{d}_{r}$ is itself a functor, sending a map $g : X \to Y$ to the map $J^{d}_{r} g : J^{d}_{r} X \to J^{d}_{r} Y$ that acts on points by post-composition. For $X$ an analytic space, there is an analogous construction that appears in \cite[\S2.3]{periodimages}. If $K \subset \mathbb{C}$ is a subfield, this construction is compatible with analytification.

The purpose of introducing jets is the following result, proven in \cite{urbanik2021sets}, building on \cite{periodimages}.

\begin{thm}
\label{jetcorresp}
For each $d, r \geq 0$, a variation of Hodge structure $\mathbb{V}$ on $S$ gives rise to a canonical map 
\[ \eta^{d}_{r} : J^{d}_{r} S \to \GL_{m} \backslash J^{d}_{r} \ch{L} , \] 
of algebraic stacks characterized by the property that for any local period map $\psi : B \to \ch{L}$ and any jet $j \in J^{d}_{r} B$ we have $\psi \circ j = \eta^{d}_{r}(j)$ modulo $\GL_{m}(\mathbb{C})$.

Moreover, if the data $(\mathcal{H}, F^{\bullet}, \nabla)$ associated to the variation $\mathbb{V}$ admits a $K$-algebraic model, the map $\eta^{d}_{r}$ is defined over $K$, and there exists an algorithm to compute the $\GL_{m}$-torsor $p^{d}_{r} : \mathcal{P}^{d}_{r} \to J^{d}_{r} S$ and the $\GL_{m}$-invariant map $\alpha^{d}_{r} : \mathcal{P}^{d}_{r} \to J^{d}_{r} \ch{L}$ which defines $\eta^{d}_{r}$ from a presentation of the data $(\mathcal{H}, F^{\bullet}, \nabla)$ in terms of finitely-presented $K$-modules.
\end{thm}

We note that the computability of the torsor $\mathcal{P}^{d}_{r}$ in \autoref{jetcorresp} has in particular the following consequence: if $\mathcal{S} \subset (\GL_{m} \backslash J^{d}_{r} \ch{L})(\mathbb{C})$ is a subset which is the image under the quotient of a constructible $L$-algebraic set $\mathcal{F} \subset J^{d}_{r} \ch{L}$, where $K \subset L$ is a computable extension, then we can compute $(\eta^{d}_{r})^{-1}(\mathcal{S})$ by computing $p^{d}_{r}((\alpha^{d}_{r})^{-1}(\mathcal{F}))$. Thus if we define \vspace{0.5em}
\begin{defn}
\label{Tconstdef}
For a constructible $L$-algebraic set $\mathcal{F} \subset J^{d}_{r} \ch{L}$, with $K \subset L$ an extension, we write
\[ \mathcal{T}^{d}_{r}(\mathcal{F}) := (\eta^{d}_{r})^{-1}(\GL_{m} \cdot \mathcal{F}) . \]
Moreover, for a type $\mathcal{C} = \mathcal{C}(W)$, we will write either $\mathcal{T}^{d}_{r}(\mathcal{C})$ or $\mathcal{T}^{d}_{r}(W)$ for the set $\mathcal{T}^{d}_{r}(J^{d}_{r} W)$.
\end{defn} \vspace{0.3em}
\noindent then the main consequence of the preceding discussion for our situation is the following, which is immediate from what we have said:
\begin{prop}
\label{diffconstprop}
For each $d, r \geq 0$ there exists an algorithm which, given a constructible $L$-algebraic set $\mathcal{F} \subset J^{d}_{r} \ch{L}$ with $K \subset L$ a computable extension, computes $\mathcal{T}^{d}_{r}(\mathcal{F}) \subset J^{d}_{r} S$. \qed
\end{prop}

\subsection{Jet Conditions and Types}

Let us now try to understand how computing the ``differential constraints'' induced by types $\mathcal{C}(W)$ as in \autoref{diffconstprop} can help us carry out Step Two of \autoref{methodsketch}. Let $\Gamma = \textrm{Aut}(V, Q)(\mathbb{Z})$, and let $\varphi : S_{\mathbb{C}} \to \Gamma \backslash D$ be the canonical period map which sends a point $s \in S(\mathbb{C})$ to the isomorphism class of the polarized Hodge structure on $\mathbb{V}_{s}$. By \cite{OMINGAGA}, the map $\varphi$ factors as $\iota \circ p$, where $p : S_{\mathbb{C}} \to T$ is a dominant map of algebraic varieties and $\iota$ is a closed embedding of analytic spaces; it follows that for each subvariety $Z \subset S_{\mathbb{C}}$ the dimension of the image $\varphi(Z)$ makes sense as the dimension of a constructible algebraic set.

Fix a sequence of compatible embeddings 
\[ \Spec K = \mathbb{D}^{0}_{r} \xhookrightarrow{\iota_{0}} \mathbb{D}^{1}_{r} \xhookrightarrow{\iota_{1}} \mathbb{D}^{2}_{r} \xhookrightarrow{\iota_{2}} \mathbb{D}^{3}_{r} \xhookrightarrow{\iota_{3}} \mathbb{D}^{4}_{r} \xhookrightarrow{\iota_{4}} \cdots \]
of formal disks. By acting on points via pullback, we obtain natural transformations of functors $\textrm{res}^{d}_{e} : J^{d}_{r} \to J^{e}_{r}$ which produce maps $J^{d}_{r} X \to J^{e}_{r} X$ that take jets $j : \mathbb{D}^{d}_{r} \to X$ to their restrictions $j \circ \iota_{d-1} \circ \cdots \circ \iota_{e}$. We are now ready to present the key proposition for our method:

\begin{defn}
For a scheme $X$ (resp. analytic space $X$) denote by $J^{d}_{r,nd} X \subset J^{d}_{r} X$ the subscheme (resp. the analytic subspace) parametrizing those maps $j : \mathbb{D}^{d}_{r} \to X$ which are injective on the level of tangent spaces. We call such $j$ \emph{non-degenerate} jets.
\end{defn}

\begin{prop}
\label{mainjetprop}
Let $\mathcal{S}$ be a set of types containing all the Hodge-theoretic types, and let $e$ and $k$ be non-negative integers. Then the following are equivalent:
\begin{itemize}
\item[(i)] there exists a geometrically irreducible subvariety $Z \subset S_{\mathbb{C}}$ with $\dim Z > d$, $\dim \varphi(Z) \geq e$, and such that $\dim \mathcal{C}(Z) - \dim \varphi(Z) \leq k$;
\item[(ii)] there exists $\mathcal{C} \in \mathcal{S}$ with $\dim \mathcal{C} - e \leq k$, and such that the intersection 
\[ \mathcal{K}^{d}_{r}(\mathcal{C}, e, k) := \mathcal{T}^{d+1}_{r}(\mathcal{C}) \cap \mathcal{T}^{d+1}_{r}((\textrm{res}^{d+1}_{e})^{-1}(J^{e}_{r,nd} \ch{L})) \cap J^{d+1}_{r,nd} S \] 
is non-empty for each $r \geq 0$.
\end{itemize}
\end{prop}

\begin{rem}
In the situation that the variation $\mathbb{V}$ admits a local Torelli theorem, one can drop the distinction between $\dim Z$ and $\dim \varphi(Z)$ and consider instead the intersections $\mathcal{T}^{d+1}_{r}(\mathcal{C}) \cap J^{d+1}_{r,nd} S$ in part (ii), ignoring the middle term.
\end{rem} \vspace{0.5em}

The rest of this section we devote to proving \autoref{mainjetprop}, identifying $S$ with $S_{\mathbb{C}}$ for ease of notation. To begin with, let us check that (i) implies (ii) by applying the definitions. If $g : S' \to S$ is an \'etale cover and we consider the variation $\mathbb{V}' = g^{*} \mathbb{V}$, then the maps $\eta^{d}_{r}$ and $\eta'^{d}_{r}$ obtained from \autoref{jetcorresp} are related by $\eta'^{d}_{r} = \eta^{d}_{r} \circ (J^{d}_{r} g)$. Choosing a finite index subgroup $\Gamma' \subset \Gamma$ and passing to such a cover, we can reduce to the case where we have a period map $\varphi : S \to \Gamma \backslash D$ with $D \to \Gamma \backslash D$ a local isomorphism. Applying \cite{OMINGAGA} the map $\varphi : S \to \Gamma \backslash D$ factors as $\varphi = \iota \circ p$, where $p : S \to T$ is a dominant map of algebraic varieties and $\iota$ is an analytic closed embedding. Then via $p$, the variety $Z$ is dominant over a closed subvariety $Y \subset T$ of dimension $\dim \varphi(Z) \geq e$. Shrinking $S$ (and hence $Z$) we may assume that $Z$ is smooth, and that $Z$ is surjective onto a dense open subset $Y^{\circ} \subset Y$. Shrinking $S$ even further we may assume that $Z \to Y^{\circ}$ is smooth. The smoothness of $Z \to Y^{\circ}$ implies in particular that the induced jet space maps $J^{d}_{r} Z \to J^{d}_{r} Y^{\circ}$ for all choices of $d$ and $r$ are surjective.

We may choose neighbourhoods $B \subset S(\mathbb{C})$ and $U \subset D$ such that $\restr{\pi}{U} : U \to \pi(U)$ is an isomorphism, both $B \cap Z$ and $\pi(U) \cap Y^{\circ}$ are non-empty, and we have a local lift $\psi : B \to U$ of $\varphi$. Choose a jet $\sigma \in J^{e}_{r,nd} (Y^{\circ} \cap \pi(U))$ and lift it along $p$ to a jet $\widetilde{\sigma} \in J^{e}_{r,nd} (Z \cap B)$ landing at the point $s \in S(\mathbb{C})$. Using the fact that the germ $(Z, s)$ is smooth of dimension $\dim Z > d$ the jet $\widetilde{\sigma}$ can be extended to a jet $j \in J^{d+1}_{r, nd} (Z \cap B)$ such that $\textrm{res}^{d+1}_{e} (j) = \widetilde{\sigma}$, and hence $\textrm{res}^{d+1}_{e} (\varphi \circ j) = \sigma$. From the fact that $\restr{\varphi}{B} = \pi \circ \psi$ and the defining property of the map $\eta^{d}_{r}$ it follows that $j$ lies inside $\mathcal{T}^{d+1}_{r}((\textrm{res}^{d+1}_{e})^{-1}(J^{e}_{r,nd} \ch{L})) \cap J^{d+1}_{r,nd} S$. We can then take $\mathcal{C} = \mathcal{C}(Z)$, and the fact that $j$ factors through $Z$ implies that $j \in \mathcal{T}^{d+1}_{r}(\mathcal{C})$ as well.

To prove the reverse implication, we review some preliminary facts relating to jets.

\begin{defn}
We say a sequence $\{ j_{r} \}_{r \geq 0}$ with $j_{r} \in J^{d}_{r} X$ is \emph{compatible} if the projections $J^{d}_{r} X \to J^{d}_{r-1} X$ map $j_{r}$ to $j_{r-1}$.
\end{defn}

\begin{lem}
\label{compseqlem}
Suppose that $\mathcal{T}_{r} \subset J^{d}_{r} X$ is a collection of non-empty constructible algebraic sets such that the projections $J^{d}_{r} X \to J^{d}_{r-1} X$ map $\mathcal{T}_{r}$ into $\mathcal{T}_{r-1}$. Then there exists a compatible sequence $\{ j_{r} \}_{r \geq 0}$ with $j_{r} \in \mathcal{T}_{r}$ for all $r \geq 0$.
\end{lem}

\begin{proof}
See \cite[Lem. 5.3]{periodimages}.
\end{proof}

\begin{defn}
Given a variety $Z$ (algebraic or analytic) and $z \in Z$ a point, we denote by $(J^{d}_{r} Z)_{z}$ the fibre above $z$ of the natural projection map $J^{d}_{r} Z \to Z$.
\end{defn}

\begin{lem}
\label{factorthrough}
If $g : (Z, z) \to (Y, y)$ is a map of analytic germs with $\dim (Z, z) = d$ and $(Z, z)$ smooth, we have an infinite compatible family $j_{r} \in (J^{d}_{r,nd} Z)_{z}$, and $g \circ j_{r} \in (J^{d}_{r} X)_{y}$ for some germ $(X, y) \subset (Y, y)$ and all $r \geq 0$, then $g$ factors through the inclusion $(X, x) \subset (Y, y)$. 
\end{lem}

\begin{proof}
See \cite[Lem. 4.5]{urbanik2021sets}.
\end{proof}

\begin{lem}
\label{jetdimlem}
Suppose that $X$ is an algebraic variety (resp. analytic space) and $x \in X$ is a point for which the fibre $(J^{d}_{r,nd} X)_{x}$ above $x$ is non-empty for all $r \geq 0$. Then the germ $(X, x)$ has dimension at least $d$. 
\end{lem}

\begin{proof}
See \cite[Prop. 2.7]{periodimages}.
\end{proof}

\begin{proof}[Proof of \ref{mainjetprop}:]
By what we have said, we are reduced to showing that (ii) implies (i). The statement is unchanged by replacing $S$ with a finite \'etale covering $g : S' \to S$ and the variation $\mathbb{V}$ with $g^{*} \mathbb{V}$; as before this does not affect the hypothesis (ii) since the maps $\eta^{d}_{r}$ and $\eta'^{d}_{r}$ associated to $S$ and $S'$ are related by $\eta'^{d}_{r} = \eta^{d}_{r} \circ (J^{d}_{r} g)$. Choosing a finite index subgroup $\Gamma' \subset \Gamma$ and choosing $g$ so the monodromy of $g^{*} \mathbb{V}$ lies in $\Gamma'$ we may reduce to the case where $D \to \Gamma \backslash D$ is a local isomorphism. Moreover, taking a futher finite \'etale cover we may apply \cite[Cor. 13.7.6]{CMS} to reduce to the case where $\varphi$ is proper; this requires possibly extending $S'$ to a variety $S''$ by adding a closed subvariety at infinity, but as long as we are careful to work only with jets that factor through $S'$ our proof will produce a variety $Z$ intersecting $S'$; in particular, we now assume that $\varphi : S \to \Gamma \backslash D$ is proper but redefine the sets $\mathcal{K}^{d}_{r}$ to equal
\[ \mathcal{T}^{d+1}_{r}(\mathcal{C}) \cap \mathcal{T}^{d+1}_{r}((\textrm{res})^{d+1}_{e})^{-1}(J^{d}_{r,nd} \ch{L})) \cap J^{d+1}_{r,nd} S^{\circ} , \]
for some open subvariety $S^{\circ} \subset S$.

Applying the main result of \cite{OMINGAGA}, the map $\varphi$ once again factors as $\varphi = \iota \circ p$ with $p : S \to T$ a dominant (now proper) map of algebraic varieties. We can then consider the Stein factorization $S \xrightarrow{q} U \xrightarrow{r} T$ of $p$; note that $q$ is proper with connected fibres, $U$ is normal, and $r$ is finite. One can define the type of a subvariety $Y \subset U$ exactly as in \autoref{vartypedef} with respect to the period map $U \to \Gamma \backslash D$. From \autoref{compseqlem} above, the assumption (ii) entitles us to a compatible sequence $\{ j_{r} \}_{r \geq 0}$ of jets such that $j_{r} \in \mathcal{K}^{d}_{r}(\mathcal{C}, e, k)$ for all $r \geq 0$. Let us write $\mathcal{C} = \mathcal{C}(W)$ for some subvariety $W \subset \ch{L}$.

By construction, the jets $\sigma_{r} = \textrm{res}^{d+1}_{e} j_{r}$ are non-degenerate, and remain so after composing with any local period map $\psi : B \to D$ for which $\sigma_{r}$ factors through $B$. This in particular implies (since $D \to \Gamma \backslash D$ is a local isomorphism) that the jets $\varphi \circ \sigma_{r}$ are non-degenerate, and hence so are the jets $q \circ \sigma_{r}$. Let $Y \subset U$ be the smallest algebraic subvariety such that $q \circ j_{r} \in J^{d+1}_{r} Y$ for all $r$. We observe that there exists a component $Z$ of $q^{-1}(Y)$ of dimension at least $d+1$ that contains the image of the jets $\{ j_{r} \}_{r \geq 0}$: one can see this by picking a neighbourhood of $j_{0}$ of the form $\mathbb{C}^{\ell} \times \mathbb{C}^{d+1}$ such that $j_{r}$ is constant on the first factor, and applying \autoref{factorthrough} above to see that the restriction of $q$ to $\{ 0 \} \times \mathbb{C}^{d+1}$ factors through $Y$. Moreover, we must have $q(Z) = Y$ by minimality, and by applying \autoref{jetdimlem} to the non-degenerate sequence $\{ q \circ \sigma_{r} \}_{r \geq 0}$ that $\dim Y \geq e$. Since $r$ is finite, this means $\dim \varphi(Z) \geq e$. From the fact that local period maps on $S$ factor through local period maps on $U$ we learn that $\mathcal{C}(Z) = \mathcal{C}(Y)$, so the result will follow if we can show that $\dim \mathcal{C}(Y) - \dim Y \leq k$. For ease of notation let us now write $\tau_{r} = q \circ j_{r}$.

Fix a local lift $\psi : B \to D$ of the period map $U \to \Gamma \backslash D$ with $B \subset U(\mathbb{C})$ an analytic ball such that the jets $\tau_{r}$ factor through $B$. Consider the set $\mathcal{G}_{r} \subset \GL_{m}(\mathbb{C})$ consisting of those $g \in \GL_{m}(\mathbb{C})$ for which $\psi \circ \tau_{r} \in g \cdot (J^{d+1}_{r} W)$. Then for each $r$ the set $\mathcal{G}_{r}$ is algebraically-constructible, and using the fact that $j_{r} \in \mathcal{T}^{d+1}_{r}(W)$ the set $\mathcal{G}_{r}$ is necessarily non-empty. Let $g_{\infty}$ be an element of this intersection. Extend $\psi$ to a lift $\widetilde{\varphi}_{Y} : \widetilde{Y} \to D$ of $Y \to \Gamma \backslash D$ to the universal covering. Then $\widetilde{\varphi}_{Y}(\widetilde{Y}) \subset D$ is a closed analytic set containing the jets $\psi \circ \tau_{r}$, and hence the non-degenerate jets $\textrm{res}^{d+1}_{e}(\psi \circ \tau_{r})$. Letting $A \subset \widetilde{\varphi}_{Y}(\widetilde{Y}) \cap (g_{\infty} \cdot W)$ be the minimal analytic germ through which $\psi \circ \tau_{r}$ (and hence $\textrm{res}^{d+1}_{e}(\psi \circ \tau_{r})$) factors, it follows from \autoref{jetdimlem} that $A$ has dimension at least $e$.

Consider the Zariski closure $V \subset Y$ of $\psi^{-1}(A)$. We claim that $V = Y$. Because $Y$ was chosen minimal containing the compatible family of jets $\{ \tau_{r} \}_{r \geq 0}$, it suffices to show that each $\tau_{r}$ factors through $V$. Consider the component of $q^{-1}(B)$ containing $j_{0}$; by choosing coordinates we may assume $q^{-1}(B) \subset \mathbb{C}^{\ell} \times \mathbb{C}^{d+1}$ is an open neighbourhood and identify $j_{0}$ with the origin. After a further change of coordinates we may assume $j_{r}$ is constant on the first factor, and let $F = q^{-1}(B) \cap (\{ 0 \} \times \mathbb{C}^{d+1})$. Applying \autoref{factorthrough} we find that $\psi(q(F)) \subset A$, and hence $q(F) \subset V$. Using proper base change the map $q^{-1}(B) \to B$ is proper, so $q(F)$ is an analytic subvariety of $B$, and by construction the jets $\{ \tau_{r} \}_{r \geq 0}$ factor through it, hence through $V$.

We are now ready to apply the Bakker-Tsimerman transcendence theorem; the jets are no longer needed. It follows from the structure theorem for period mappings \cite[III.A]{GGK}, the closed analytic set $\widetilde{\varphi}_{Y}(\widetilde{Y})$ lies inside an orbit $\ch{D}' = \mathbf{H}_{Y} \cdot \psi(\tau_{0})$ of the algebraic monodromy of $Y$. Consider the graph $E$ of $\widetilde{\varphi}_{Y}$ in $Y \times \ch{D}'$. Then as $A$ has dimension $e$ and $\psi^{-1}(A)$ is Zariski dense, there exists a component $C$ of $E \cap (Y \times (\ch{D}' \cap g_{\infty} \cdot W))$ of dimension at least $e$ and projecting to a Zariski dense subset of $Y$. Applying the main theorem of \cite{AXSCHAN} we learn that
\begin{align*}
\textrm{codim}_{Y \times \ch{D}'} (Y \times (\ch{D}' \cap g_{\infty} \cdot W)) + \textrm{codim}_{Y \times \ch{D}'} E &\leq \textrm{codim}_{Y
 \times \ch{D}'} C \\
(\dim \ch{D}' - \dim W) + \dim \ch{D}' &\leq \dim Y + \dim \ch{D}' - \dim C \\
\dim \ch{D}' - \dim Y &\leq \dim W - \dim C \\
\dim \mathcal{C}(Y) - \dim Y &\leq \dim \mathcal{C} - e \\
&\leq k
\end{align*}
\end{proof}

\section{Main Results}

\subsection{Computing Bounds on $\Delta_{d}$}

\subsubsection{Computing Lower Bounds}

Let us explain the significance of \autoref{mainjetprop} in proving \autoref{mainthm}, i.e., giving an effective method to compute bounds for \[ \Delta_{d} = \min_{\dim Z > d} [ \dim \mathcal{C}(Z) - \dim \varphi(Z) ] , \]
where we have used \autoref{mainthm2} and \autoref{vartypedef} to give this equivalent expression for $\Delta_{d}$. Since $\Delta_{d}$ is a integer bounded by $\dim D$, giving an effective method to compute it amounts to developing a procedure to decide, for any integer $k$, whether we have $\Delta_{d} \leq k$. This in turn amounts to deciding, for any integer $0 \leq e \leq \dim \varphi(S)$, whether (ii) holds in \autoref{mainjetprop}.

Let us take $\mathcal{S}$ to be the set up types computed by Step One, and let us suppose that in fact $\Delta_{d} > k$. Then by the equivalence in \autoref{mainjetprop}, we should find that for any $\mathcal{C} \in \mathcal{S}$ with $\dim \mathcal{C} - e \leq k$, there must be some $r = r(\mathcal{C}, e, k)$ such that $\mathcal{K}^{d}_{r}(\mathcal{C}, e, k)$ is empty. Moreover, verifying that such an $r$ exists for each such $\mathcal{C}$ and $e$ proves, again by the same equivalence, that $\Delta_{d} > k$. Consequently, we obtain the following result, which is the first half of \autoref{mainthm}:

\begin{prop}
\label{lowerboundcomp}
By computing the sets $\mathcal{K}^{d}_{r}(\mathcal{C}, e, k)$ described in \autoref{mainjetprop} in parallel, we may compute a non-decreasing sequence of lower bounds
\[ \kappa(1) \leq \kappa(2) \leq \cdots \leq \kappa(r) \leq \cdots < \Delta_{d} \]
such that for some $r = r_{0}$ we have $\kappa(r_{0}) = \Delta_{d} - 1$.
\end{prop}

\begin{proof}
At the $r$'th stage we compute all the sets $\mathcal{K}^{d}_{r}(\mathcal{C}, e, k)$ for all applicable choices of $\mathcal{C}$, $e$ and $k$, and set $\kappa(r)$ to be the smallest $k$ for which all the sets $\mathcal{K}^{d}_{r}(\mathcal{C}, e, k)$ are empty. From the discussion preceding the Proposition, the result follows.
\end{proof}

\subsubsection{Computing Upper Bounds}

\autoref{lowerboundcomp} does not actually give an algorithm for computing $\Delta_{d}$, since no way is given to decide when $r = r_{0}$. For applications to the Lawrence-Venktesh method this doesn't matter: one wants to be able to compute an optimal lower bound for $\Delta_{d}$, but one does not actually have to prove that this lower bound actually equals $\Delta_{d}$ in order to apply the diophantine finiteness machinery. Nevertheless, let us explain how one can do this in the case where $\varphi$ is quasi-finite; under this assumption, we may drop the distinction between $\dim Z$ and $\dim \varphi(Z)$, and we are instead interested in computing
\[ \min_{\dim Z > d} \, \left[ \dim \mathcal{C}(Z) - \dim Z \right] . \]

 What is needed is the following:

\begin{prop}
\label{upperboundcomp}
Suppose that $S$ is quasi-projective and $\varphi$ is quasi-finite. Then there exists a procedure that outputs an infinite sequence of upper bounds
\[ \tau(1) \geq \tau(2) \geq \cdots \geq \tau(i) \geq \cdots \geq \Delta_{d} \] 
such that for some $i = i_{0}$ we have $\tau(i) = \Delta_{d}$.
\end{prop}

\noindent Given both \autoref{lowerboundcomp} and \autoref{upperboundcomp} we obtain an algorithm for computing $\Delta_{d}$ by running both procedures in parallel and terminating when $\kappa(r) + 1 = \tau(i)$. 

\subsection{Finding Varieties that Exhibit $\Delta_{d}$}

In this section we prove \autoref{upperboundcomp}, assuming throughout that $S$ is quasi-projective and $\varphi$ is quasi-finite. Let us fix a projective compactification $S \subset \overline{S}$ of $S$ and consider the Hilbert scheme $\textrm{Hilb}(\overline{S})$. There exist algorithms, for instance by working with the Pl\"uker coordinates of the appropriate Grassmannian, for computing any finite subset of components of $\textrm{Hilb}(\overline{S})$. By \cite[Lem. 5.10]{urbanik2021sets} we obtain the same fact for the open locus $\textrm{Var}(S) \subset \textrm{Hilb}(\overline{S})$ consisting of just those points $[\overline{Z}]$ for which $Z = S \cap \overline{Z}$ is a non-empty geometrically irreducible algebraic subvariety of $S$. What we will show is that there exists a procedure which outputs an infinite sequence $\{ \mathcal{W}_{i} \}_{i = 1}^{\infty}$ of constructible algebraic loci $\mathcal{W}_{i} \subset \textrm{Var}(S)$, with the following two properties:
\begin{itemize}
\item[(i)] for each $i$, the type $\mathcal{C}(Z)$ and dimension $\dim Z$ are constant over all $[Z] \in \mathcal{W}_{i}$;
\item[(ii)] there exists some $i = i_{0}$ such that 
\[ \Delta_{d} = \dim \mathcal{C}(Z) - Z , \]
for some (hence any) point $[Z] \in \mathcal{W}_{i}$.
\end{itemize}
Given such an algorithm the problem of computing the bound $\tau(i)$ that appears in \autoref{upperboundcomp} reduces to choosing a point $[Z] \in \mathcal{W}_{i}$, computing $\dim \mathcal{C}(Z) - \dim Z$, and setting 
\[ \tau(i) := \textrm{min} \{ \tau(i-1), \dim \mathcal{C}(Z) - \dim Z \} . \] 
(We note that the problem of computing $\dim \mathcal{C}(Z)$ from $Z$ and the restriction $\restr{(\mathcal{H}, F^{\bullet}, \nabla)}{Z}$ of the algebraic data on $S$ is solved for us by \cite[Lem. 5.8]{urbanik2021sets} by taking the family $g$ in the statement of \cite[Lem. 5.8]{urbanik2021sets} to be a trivial family; we will say little about this problem here.)

In fact, an algorithm for computing the sets $\mathcal{W}_{i}$ has already been given in a previous paper by the author. We begin by recalling the necessary background. We regard $S$ as a complex algebraic variety in what follows. Given two (geometrically) irreducible subvarieties $Z_{1}, Z_{2} \subset S$ with $Z_{1} \subset Z_{2}$, the algebraic monodromy group $\mathbf{H}_{Z_{1}}$ may be naturally regarded as a subgroup of $\mathbf{H}_{Z_{2}}$ (after choosing a base point $s \in Z_{1}(\mathbb{C})$). Using this we define:
\begin{defn}
An irreducible complex subvariety $Z \subset S$ is said to be \emph{weakly special} if it is maximal among such subvarieties for its algebraic monodromy group.
\end{defn}
The key fact is then the following:
\begin{lem}
For each integer $d > 0$, there exists a weakly special subvariety $Z \subset S$ such that 
\[ \Delta_{d} = \dim \mathcal{C}(Z) - \dim Z . \]
\end{lem}

\begin{proof}
By \cite[Prop. 4.18]{urbanik2021sets}, the condition that $Z$ be weakly special is equivalent to $Z$ being a maximal irreducible complex subvariety of $S$ of type $\mathcal{C}(Z)$. Thus if we have any $Y$ which is not weakly special, there exists a weakly special $Z$ properly containing $Y$ with $\mathcal{C}(Y) = \mathcal{C}(Z)$, hence
\[ \dim \mathcal{C}(Z) - \dim Z < \dim \mathcal{C}(Y) - \dim Y . \]
It follows that the value of $\Delta_{d}$ can only be achieved by a weakly special variety.
\end{proof}

\begin{proof}[Proof of \ref{upperboundcomp}]
By the \emph{degree} of a subvariety $Z \subset S$ we will mean the degree of its closure $\overline{Z}$ inside $\overline{S}$. For any integer $b$, denote by $\textrm{Var}(S)_{b} \subset \textrm{Var}(S)$ the finite-type subscheme parametrizing varieties of degree at most $b$. Denote by $\mathcal{W} \subset \textrm{Var}(S)$ the locus of weakly special subvarieties. Then given an integer $b$, the algorithm that appears in \cite[Thm. 5.15]{urbanik2021sets} computes the intersection $\mathcal{W} \cap \textrm{Var}(S)_{b}$ as a constructible algebraic locus.

Let us describe the algorithm appearing in \cite[Thm. 5.15]{urbanik2021sets} more precisely. Consider the types $\mathcal{C}_{1}, \hdots, \mathcal{C}_{\ell}$ computed by Step One, and define for each such type $\mathcal{C}_{j}$ the locus
\[ \mathcal{W}(\mathcal{C}_{j}) := \{ [Z] \in \textrm{Var}(S) : \mathcal{C}(Z) \leq \mathcal{C}_{j} \} . \] 
It is shown in \cite[Prop. 4.31]{urbanik2021sets} that for each $j$ the locus $\mathcal{W}(\mathcal{C}_{j})$ is closed algebraic. We can then consider the sublocus $\mathcal{W}(\mathcal{C}_{j})_{\textrm{opt}} \subset \mathcal{W}(\mathcal{C}_{j})$ consisting of just those components $C \subset \mathcal{W}(\mathcal{C}_{j})$ for which a generic point $[Z] \in C$ satisfies $\mathcal{C}(Z) = \mathcal{C}_{j}$.

In \cite[Prop. 5.14]{urbanik2021sets}, an algorithm is given for computing $\mathcal{W}(\mathcal{C}_{j})_{\textrm{opt}} \cap \textrm{Var}(S)_{b}$ for each $j$. Using this, one can compute all the finitely many closed algebraic loci $C_{1}, \hdots, C_{i_{b}}$ which arise as a component of $\mathcal{W}(\mathcal{C}_{j})_{\textrm{opt}} \cap \textrm{Var}(S)_{b}$ for some $j$. The problem of computing $\mathcal{W} \cap \textrm{Var}(S)_{b}$ is reduced to computing constructible algebraic conditions on each component $C_{i} \subset \mathcal{W}(\mathcal{C}_{j})_{\textrm{opt}}$ which define the locus $\mathcal{W}_{i} \subset C_{i}$ of points $[Z] \in C_{i}$ that are weakly special of type $\mathcal{C}_{j}$. This is taken care of by the proof of \cite[Thm. 5.15]{urbanik2021sets}. By construction, the points in $\mathcal{W}_{i}$ all have the same type and the same dimension, so we complete the proof by computing these loci for increasing values of $b$.
\end{proof}

\section{Application to Lawrence-Venkatesh}

We now show how the bound of \autoref{mainthm} can be used to establish diophantine finiteness results. Similar arguments appear in \cite{LV} and \cite{lawrence2020shafarevich}, but as they are not precisely adapted to our setup, we give our own version. We recall the situation: we have a smooth projective family $f : X \to S$ over the smooth base $S$, with everything defined over $\mathcal{O}_{K,N}$.\footnote{Note in particular we are assuming now that $S$ is smooth over $\mathcal{O}_{K,N}$, which we can achieve by increasing $N$ if necessary.} The relative algebraic de Rham cohomology $\mathcal{H} = R^{i} f_{*} \Omega^{\bullet}_{X/S}$ gives a model for the Hodge bundle $\mathbb{V} \otimes \mathcal{O}_{\an{S}}$, where $\mathbb{V} = R^{i} f_{*} \mathbb{Z}$. By a result of Katz and Oda \cite{katz1968}, the flat connection associated to the local system $\mathbb{V}_{\mathbb{C}}$ by the Riemann-Hilbert correspondence admits a model $\nabla : \mathcal{H} \to \Omega^{1}_{S} \otimes \mathcal{H}$ after possibly increasing $N$. Likewise, we may also assume the Hodge filtration $F^{\bullet}$ gives a filtration of $\mathcal{H}$ by vector subbundles. 

Fix a prime $p$ not dividing $N$, and a place $v$ of $K$ above $p$. Then for each integral point $s \in S(\mathcal{O}_{K,N})$, we have a Galois representation $\rho_{s} : \textrm{Gal}(\overline{K}/K) \to \textrm{Aut}(H^{i}_{\textrm{\'et}}(X_{\overline{K}, s}, \mathbb{Q}_{p}))$, and an argument of Faltings \cite[Lem 2.3]{LV} shows that the semisimplifications of the representations $\rho_{s}$ belong to a finite set of isomorphism classes. From crystalline cohomology, each $s \in S(\mathcal{O}_{K,N})$, viewed as a point of $S(\mathcal{O}_{K,v})$ where $\mathcal{O}_{K,v}$ is the $v$-adic ring of integers, gives rise to a triple $(H^{i}_{\textrm{dR}}(X_{s}), \phi_{s}, F^{\bullet}_{s})$ where $\phi_{s}$ is the crystalline Frobenius. Moreover, using the functor $D_{\textrm{cris}}$ of $p$-adic Hodge theory \cite[Expos\'e III]{fontaine1994corps}, the triple $(H^{i}_{\textrm{dR}}(X_{s}), \phi_{s}, F^{\bullet}_{s})$ is determined up to isomorphism by the restriction $\rho_{s,v}$ along the map $\textrm{Gal}(\overline{K_{v}}/K_{v}) \to \textrm{Gal}(\overline{K}/K)$ determined by a fixed embedding $\overline{K} \hookrightarrow \overline{K_{v}}$. We denote by $\mathcal{I}(s)$ all those triples $(V, \phi, F^{\bullet})$ which are of the form $D_{\textrm{cris}}(\restr{\rho}{\mathbb{Q}_{p}})$, where $\rho$ is a global Galois representation whose semisimplification is isomorphic to the semisimplification of $\rho_{s}$.

Recall that we have fixed the integral lattice $V = \mathbb{Z}^{m}$, where $m$ is the dimension of the cohomology of the fibres of $f$, and a $\mathbb{Q}$-algebraic flag variety $\ch{L}$ of Hodge flags on $V$. In what follows we write $V_{p}$ for $V \otimes \mathbb{Q}_{p}$, and $\ch{L}_{p}$ for $\ch{L}_{\mathbb{Q}_{p}}$. Then the key idea of the Lawrence-Venkatesh method is the following:

\begin{prop}
\label{LVprop}
Suppose that for each $s \in S(\mathcal{O}_{K,N})$, whenever we have an endomorphism $\phi_{s} : V_{p} \to V_{p}$ and a flag $F^{\bullet}_{s}$ on $V_{p}$ such that $(V_{p}, \phi_{s}, F^{\bullet}_{s})$ represents $\mathcal{I}(s)$, the Hodge flags $F^{\bullet}$ on $V_{p}$ for which $(V_{p}, \phi_{s}, F^{\bullet}) \in \mathcal{I}(s)$ lie in an algebraic subvariety $O_{s} \subset \ch{L}_{p}$ satisfying $\Delta_{d} \geq \dim O_{s}$. Then $\dim \overline{S(\mathcal{O}_{K,N})}^{\textrm{Zar}} \leq d$.  
\end{prop}

To prove \autoref{LVprop} we will need a rigid-analytic version of the Bakker-Tsimerman transcendence theorem, which we will see can be deduced formally from the complex analytic one. To set things up, let us revisit the term \emph{local period map}, this time in the rigid analytic setting (c.f. \autoref{locperdef}). We will denote by $\mathbb{C}_{p}$ the completion of the algebraic closure $\overline{K_{v}}$. In what follows we sometimes identify algebraic varieties with their rigid-analytifications when the context is clear.

\begin{defn}
\label{padiclocperdef}
Let $K_{p}$ be a local field containing $K_{v}$, let $\an{S_{K_{p}}}$ be the rigid-analytification of the base-change $S_{K_{p}}$ of $S$, and suppose that $B_{p} \subset \an{S_{K_{p}}}$ is an affinoid subdomain. Then a (rigid-analytic) local period map $\psi : B_{p} \to \an{\ch{L}_{K_{p}}}$ is a rigid-analytic map obtained as a composition $\psi = \an{q_{K_{p}}} \circ A_{p}$, where:
\begin{itemize}
\item[(i)] The rigid analytifications $F^{k} \an{\mathcal{H}_{K_{p}}}$ are all trivial on $B_{p}$.
\item[(ii)] The map $A_{p} : B_{p} \to \an{\GL_{m, K_{p}}}$ is a varying filtration-compatible $p$-adic period matrix over $B_{p}$. More precisely, there exists a basis $v^{1}, \hdots, v^{m}$ for $\an{\mathcal{H}_{K_{p}}}(B_{p})$, compatible with the filtration in the sense that $F^{k} \an{\mathcal{H}_{K_{p}}}(B_{p})$ is spanned by $v^{1}, \hdots, v^{i_{k}}$ for some $i_{k}$, and a flat (for $\an{\nabla_{K_{p}}}$) frame $b^{1}, \hdots, b^{m}$ such that $A_{p}$ gives a varying change-of-basis matrix from $v^{1}, \hdots, v^{m}$ to $b^{1}, \hdots, b^{m}$.
\item[(iii)] The map $q : \GL_{m} \to \ch{L}$ is the map that sends a matrix $M$ to the Hodge flag $F^{\bullet}_{M}$ defined by the property that $F^{k}_{M}$ is spanned by the first $i_{k}$ columns.
\end{itemize}
\end{defn}

To prove \autoref{LVprop} we will need a version of the Bakker-Tsimerman transcendence result for rigid-analytic local period maps, which we prove by formally transferring the same result for complex analytic local period maps. To avoid certain minor pathologies that can occur in the complex analytic case we will restrict to local period maps $\psi : B \to \ch{L}$ which are definable in the structure $\mathbb{R}_{\textrm{an}, \textrm{exp}}$; for background on definability and definable analytic spaces we refer to \cite{van1996geometric} and \cite{OMINGAGA}. We note that this is not a serious restriction: given any local period map $\psi$ and any point $s \in B$ there exists a definable restriction of $\psi$ to a neighbourhood of $s$, a fact which is for instance easily deduced from \cite[Prop. 4.27]{urbanik2021sets}.

\begin{lem}
\label{padicaxschanlem}
~\begin{itemize}
\item[(i)] Suppose that $\psi : B \to \an{\ch{L}_{\mathbb{C}}}$ is a definable analytic local period map on $\an{S_{\mathbb{C}}}$. Let $V \subset \ch{L}_{\mathbb{C}}$ be an algebraic subvariety satisfying $\Delta_{d} \geq \dim V$. Then $\psi^{-1}(V)$ lies in an algebraic subvariety of $S_{\mathbb{C}}$ of dimension at most $d$.
\item[(ii)] Suppose that $\psi_{p} : B_{p} \to \an{\ch{L}_{\mathbb{C}_{p}}}$ is a rigid-analytic local period map on $\an{S_{\mathbb{C}_{p}}}$. Let $V_{p} \subset \ch{L}_{p, \mathbb{C}_{p}}$ be an algebraic subvariety satisfying $\Delta_{d} \geq \dim V_{p}$. Then $\psi_{p}^{-1}(V_{p})$ lies in an algebraic subvariety of $S_{\mathbb{C}_{p}}$ of dimension at most $d$.
\end{itemize}
\end{lem}

\paragraph{Proof of \autoref{padicaxschanlem}(i):} ~ \\

\vspace{-0.5em}

This is an application of the Bakker-Tsimerman transcendence theorem. Let $Z \subset S_{\mathbb{C}}$ be the Zariski closure of $\psi^{-1}(V)$. We assume for contradiction that $\dim Z > d$, and let $Z_{0} \subset Z$ be a component of maximal dimension. Let $\varphi : \an{S_{\mathbb{C}}} \to \Gamma \backslash D$ be the canonical period map with $\Gamma = \textrm{Aut}(V,Q)(\mathbb{Z})$. The statement is invariant under replacing $\psi$ with a $\GL_{m}(\mathbb{C})$-translate $g \cdot \psi$ and $V$ with $g \cdot V$, so we may assume that $\psi$ is a local lift of $\varphi$. Arguing as in \cite[Cor. 13.7.6]{CMS} we may assume that $\varphi$ is proper, hence the image $T = \varphi(S)$ is algebraic by \cite{OMINGAGA}, and we may consider the Stein factorization $S \xrightarrow{q} U \xrightarrow{r} T$ of the map $S \to T$. 

Let $Y = q(Z_{0})$, and note that $\dim Y = \dim \varphi(Z_{0})$. By assumption we have $\Delta_{d} \leq \dim \mathcal{C}(Z_{0}) - \dim \varphi(Z_{0}) = \dim \mathcal{C}(Y) - \dim Y$, where the type of $Y$ is taken with respect to the period map $U \to \Gamma \backslash D$. Moreover, this continues to hold if we replace $Y$ with a smooth resolution $Y'$. The variation of Hodge structure on $S$ descends to $U$, and hence shrinking $B$ if necessary we may factor $\psi$ through a definable local lift on $U$. By pulling back along the resolution we obtain a definable local lift $\psi' : B' \to D$ of the period map $\varphi' : Y' \to \Gamma \backslash D$ such that $\psi'^{-1}(V)$ is Zariski dense in $Y'$. We are reduced to the following situation: we have smooth variety $Y'$ with a period map $\varphi' : Y' \to \Gamma \backslash D$, a local lift $\psi' : B' \to D$ such that $\psi'^{-1}(V)$ is Zariski dense, and such that $\dim \mathcal{C}(Y') - \dim Y' \geq \dim V$. 

We now contradict the Bakker-Tsimerman theorem. In particular, we may extend the local lift $\psi'$ to a lift $\widetilde{\varphi'} : \widetilde{Y'} \to D$ of $\varphi'$ to the universal cover, and consider the graph $W \subset Y' \times \ch{D}'$ of the map $\widetilde{\varphi'}$, where $\ch{D}'$ is the orbit $\mathbf{H}_{Y'} \cdot \psi'(y)$ for some $y \in Y'(\mathbb{C})$. We then have that
\begin{align*}
\textrm{codim}_{Y' \times \ch{D}'} (Y' \times (V \cap \ch{D}')) + \textrm{codim}_{Y' \times \ch{D}'} W &\geq \dim \ch{D}' - \dim V + \dim \ch{D}' \\
&= \dim \mathcal{C}(Y') - \dim V + \dim \mathcal{C}(Y') \\
&\geq \dim Y' + \dim \mathcal{C}(Y') .
\end{align*}
Since $\psi'$ is definable, $\psi'^{-1}(V)$ has finitely many components, and hence there exists an analytic component $C \subset B'$ of $\psi'^{-1}(V)$ such that $C$ is Zariski dense in $Y'$. Let $\widetilde{C} \subset Y' \times \ch{D}'$ be its graph under $\psi'$. If $\dim Y' = 0$ there is nothing to show, so we may assume that $\dim \mathcal{C}(Y') > 0$. Hence we find that $\dim Y' + \dim \mathcal{C}(Y') > \codim_{Y' \times \ch{D}'} \widetilde{C}$, and by the Bakker-Tsimerman theorem \cite{AXSCHAN} the component $C$ lies in a proper subvariety of $Y'$, giving a contradiction. \qed

\vspace{1em}

To prove \autoref{padicaxschanlem}(ii) we first translate \autoref{padicaxschanlem}(i) into a claim about rings of formal power series. In particular let $\psi : B \to \ch{L}$ be a local period map with $V \subset \ch{L}$ an algebraic subvariety, and choose a point $s \in B$ such that $\psi(s) = t \in V$. Then $\psi$ induces a map on formal power series rings $\widehat{\psi}^{\sharp} : \widehat{\mathcal{O}}_{\ch{L}_{\mathbb{C}}, t} \to \widehat{\mathcal{O}}_{S_{\mathbb{C}}, s}$. The claim of \autoref{padicaxschanlem}(i) then says that if $I_{V} \subset \mathcal{O}_{\ch{L}_{\mathbb{C}}, t}$ is the ideal defining $V$ with extension $\widehat{I}_{V}$ inside $\widehat{\mathcal{O}}_{\ch{L}_{\mathbb{C}}, t}$, then the ideal generated by $\widehat{\psi}^{\sharp}(\widehat{I}_{V})$ contains an ideal $\widehat{I}_{Z}$ which is the extension of an ideal $I_{Z} \subset \mathcal{O}_{S_{\mathbb{C}}, s}$ defining the germ of a subvariety of dimension at most $d$.

\paragraph{Proof of \autoref{padicaxschanlem}(ii):} ~ \\

\vspace{-0.5em}

The claim is Zariski-local on $S$, so we can in particular assume that the bundles $F^{k} \mathcal{H}$ for varying $k$ are algebraically trivial over $S$, that $S$ is affine, and by smoothness that $\Omega^{1}_{S}$ is free. By definition, the map $\psi : B_{p} \to \an{\ch{L}_{\mathbb{C}_{p}}}$ is associated to the following data: a filtration-compatible frame $v^{1}, \hdots, v^{m}$, where $v^{1}, \hdots, v^{i_{k}}$ spans $F^{k} \an{\mathcal{H}_{\mathbb{C}_{p}}} (B_{p})$, and a flat frame $b^{1}, \hdots, b^{m}$ spanning $\an{\mathcal{H}_{\mathbb{C}_{p}}}(B_{p})$, where flatness means $\an{\nabla_{\mathbb{C}_{p}}} b_{i} = 0$ for all $1 \leq i \leq m$. This data satisfies the property that $\psi = \an{q_{\mathbb{C}_{p}}} \circ A_{p}$, where $A_{p}$ is the change-of-basis matrix from the frame  $v^{1}, \hdots, v^{m}$ to the frame $b^{1}, \hdots, b^{m}$, and $q$ is the map $\an{\GL_{m,\mathbb{C}_{p}}} \to \an{\ch{L}_{\mathbb{C}_{p}}}$ sending a matrix to the Hodge flag it represents. We note that changing the frame $v^{1}, \hdots, v^{m}$ to another filtration-compatible frame $v'^{1}, \hdots, v'^{m}$ does not change the local period map: such a change has the effect of replacing the map $A_{p}$ with $A_{p} \cdot C$, where $C$ is a varying matrix over $B_{p}$ whose right-action on $A_{p}$ preserves the span of the first $i_{k}$ columns for each $k$, and hence $q \circ (A_{p} \cdot C) = q \circ A_{p}$. We threfore lose no generality by assuming the filtration-compatible frame is the restriction to $B_{p}$ of an algebraic filtration-compatible frame over $S$. 

The affinoid neighbourhood $B_{p}$ is of the form $\textrm{Sp} \, T$, where $T$ is an affinoid $\mathbb{C}_{p}$-algebra. The inverse image $\psi^{-1}(V)$ is then a closed affinoid subdomain of $B_{p}$, i.e., it corresponds to an ideal $I \subset T$ such that $\psi^{-1}(V)$ may be identified with $\textrm{Sp} \, T/I$. If $R$ is the coordinate ring of $S_{\mathbb{C}_{p}}$, then the map $B_{p} \hookrightarrow \an{S_{\mathbb{C}_{p}}} \to S_{\mathbb{C}_{p}}$ induces a map $\iota : R \to T$, and the claim to be shown is that there exists an ideal $J \subset R$ defining a subvariety of dimension at most $d$ such that $\iota(J) \subset I$. The ring $T$ is Noetherian, so the ideal $I$ admits a primary decomposition. Taking radicals, we obtain finitely many prime ideals $I_{1}, \hdots, I_{\ell}$ containing $I$ such that the problem reduces, for each $1 \leq j \leq \ell$, to finding $J_{j} \subset R$ defining subvarieties of dimension at most $d$ such that $\iota(J_{j}) \subset I_{j}$ for each $j$. The analytification map $\an{S_{\mathbb{C}_{p}}} \to S_{\mathbb{C}_{p}}$ is bijective onto the set of closed points of $S_{\mathbb{C}_{p}}$ and induces isomorphisms on completed local rings [see whatever]. It follows that if we choose a maximal ideal $\mathfrak{m} \subset T$ containing $I_{j}$ we obtain a commuting diagram
\begin{center}
\begin{tikzcd}
R \arrow[r,"\iota"] \arrow[d, hook] & \arrow[d, hook] T \\
\widehat{R}_{\iota^{-1}(\mathfrak{m})} \arrow[r,"\sim"] & \widehat{\mathcal{O}}_{B_{p}, \mathfrak{m}} ,
\end{tikzcd}
\end{center}
where the bottom arrow is an isomorphism of completed local rings, and the vertical arrows are injections. In particular, if we denote by $\widehat{I}_{j}$ the extension of $I_{j}$ in $\widehat{\mathcal{O}}_{B_{p}, \mathfrak{m}}$, it suffices to show that $\widehat{\iota}(J_{j}) \subset \widehat{I}_{j}$, where $\widehat{\iota}$ is the composition of the left and bottom arrow; here we have used the fact that $I_{j} = \widehat{I}_{j} \cap T$. 

Fix an isomorphism $\tau : \mathbb{C}_{p} \xrightarrow{\sim} \mathbb{C}$, which we choose to preserve the embeddings $K \subset \mathbb{C}$ and $K \subset \mathbb{C}_{p}$. Using the model for $S$ over $K$, the isomorphism $\tau$ allows us to identify $R$ with the coordinate ring of $S_{\mathbb{C}}$, the ideal $\iota^{-1}(\mathfrak{m})$ with a complex point $s \in S(\mathbb{C})$, the ring $\widehat{R}_{\iota^{-1}(\mathfrak{m})}$ with the completed local ring $\widehat{\mathcal{O}}_{S_{\mathbb{C}}, s}$. Let $t_{p}$ be the image of the point corresponding to $\mathfrak{m}$ under $\psi$, and let $t$ be the composition $t_{p} \circ \tau^{-1}$. Applying the isomorphism $\tau$ at the level of formal power series, the rigid-analytic local period map $\psi$ induces a map 
\[ \widehat{\mathcal{O}}_{\ch{L}_{\mathbb{C}}, t} \xrightarrow{\tau} \widehat{\mathcal{O}}_{\ch{L}_{\mathbb{C}_{p}}, t_{p}} \xrightarrow{\widehat{\psi}^{\sharp}} \widehat{\mathcal{O}}_{B_{p}, \mathfrak{m}} \xrightarrow{\sim} \widehat{R}_{\iota^{-1}(\mathfrak{m})} \xrightarrow{\tau} \widehat{\mathcal{O}}_{S_{\mathbb{C}}, s} , \]
whose composition we denote by $\widehat{\eta}$. In what follows we identify the ideals $\widehat{I}_{j}$ with their images in $\widehat{\mathcal{O}}_{S_{\mathbb{C}}, s}$; by construction they are the extensions along $\widehat{\eta}$ of an ideal in $\widehat{\mathcal{O}}_{\ch{L}_{\mathbb{C}}, t}$ associated to the base-change of $V$ using $\tau$. By part (i) of this theorem and our reformulation of it in terms of completed local rings, it suffices to show that $\widehat{\eta}$ is induced by a complex analytic local period map defined on a neighbourhood of $s$.

Recall that we have a decomposition $\psi = \an{q_{\mathbb{C}_{p}}} \circ A_{p}$, where $q$ is the rigid-analytification of a $\mathbb{Q}$-algebraic map, and $A_{p}$ gives a varying change-of-basis matrix between a filtration-compatible frame $v^{1}, \hdots, v^{m}$ and a rigid-analytic flat frame. Recall also that we have chosen $v^{1}, \hdots, v^{m}$ so that it is the rigid-analytification of a $K$-algebraic filtration-compatible frame $w^{1}, \hdots, w^{m}$ over $S$. Using the decomposition $\psi = \an{q_{\mathbb{C}_{p}}} \circ A_{p}$ and the isomorphism $\tau$ we may factor $\widehat{\eta}$ as $\widehat{q} \circ \widehat{\kappa}$, where $\widehat{\kappa} : \widehat{\mathcal{O}}_{S_{\mathbb{C}}, s} \to \widehat{\mathcal{O}}_{\GL_{m, \mathbb{C}}, P}$ is the base-change under $\tau$ of the map induced by $A_{p}$. From our definition of local period map in \autoref{locperdef}, it suffices to show that $\widehat{\kappa}$ is induced by a varying change-of-basis matrix $A : B \to \an{\GL_{m,\mathbb{C}}}$ from $w^{1}, \hdots, w^{m}$ to a complex-analytic flat frame.

The result will follow from the fact that $A$ and $A_{p}$ satisfy a common set of $K$-algebraic differential equations whose solutions are uniquely determined by the period matrix they assign to a point in $B$. To see this, let us write $\nabla w^{i} = \sum_{j = 1}^{m} c_{ij} \otimes w^{j}$ for $K$-algebraic sections $c_{ij} \in \Omega^{1}_{S_{K}}$. Suppose then that $b^{k} = \sum_{i = 1}^{m} f_{ik} w^{i}$ is a flat frame on some complex analytic or rigid-analytic neighbourhood. We then have that
\begin{align*}
\nabla b^{k} &= \nabla \left( \sum_{i = 1}^{m} f_{ik} w^{i} \right) \\
&= \sum_{j = 1}^{m} df_{jk} \otimes w^{j} + \sum_{i = 1}^{m} f_{ik} \left( \sum_{j = 1}^{m} c_{ij} \otimes v^{j} \right) \\
&= \sum_{j = 1}^{m} \left( df_{jk} + \sum_{i = 1}^{m} f_{ik} c_{ij} \right) \otimes v^{j} ,
\end{align*}
from which we see that $b^{k}$ giving a flat frame is equivalent to $f_{jk}$ satisfying the system of differential equations $df_{jk} = -\sum_{i = 1}^{m} f_{ik} c_{ij}$ for all $1 \leq j, k \leq m$. If we choose a trivialization $dz_{1}, \hdots, dz_{n}$ of $\Omega^{1}_{S_{K}}$, we may write the $c_{ij}$ in terms of their coefficients $c_{ij,\ell}$ with respect to this trivialization, and the same system of differential equations becomes 
\begin{equation}
\label{diffeq}
\partial_{\ell} f_{jk} = -\sum_{i} f_{ik} c_{ij,\ell} ;
\end{equation} 
here the operator $\partial_{\ell}$ is defined using the dual basis to $dz_{1}, \hdots, dz_{n}$. By differentiating \autoref{diffeq} and substituting the lower-order differential equations into the higher-order ones, we obtain, for each sequence $\{ \ell_{i} \}_{i = 1}^{e}$ with $1 \leq \ell_{i} \leq n$ and $e \geq 1$, a set of $K$-algebraic polynomials $\xi_{\ell_{1}, \hdots, \ell_{e}; jk}$ in the functions $f_{uv}$ for $1 \leq u, v \leq m$ with coefficients in the coordinate ring of $S_{K}$ such that 
\begin{equation}
\label{diffeqpolys}
\partial_{\ell_{1}} \cdots \partial_{\ell_{e}} f_{jk} = \xi_{\ell_{1}, \hdots, \ell_{e}; jk}([f_{uv}]) .
\end{equation}

Because $S_{K}$ is smooth, given a point $s$ of $S_{K}$ the functions $z_{1} - s_{1}, \hdots, z_{n} - s_{n}$, where $s_{i}$ is the value of $z_{i}$ on $s$, induce a coordinate system in the local and formal power series rings associated to $S_{K}$ at $s$. In these coordinates, the map $A_{p}$ is given by $f^{-1}_{p}$, where $f_{p} = [f_{uv}]$ is a rigid-analytic matrix-valued solution to the differential equations \autoref{diffeq}. The formal map $\widehat{\kappa}$ obtained using the isomorphism $\tau$ then satisfies the same set of differential equations, and in particular its derivatives of all orders at $s$ are determined using \autoref{diffeqpolys} by the initial condition $f^{-1}(s) = P$. If we then construct an analytic solution to the differential  system in \autoref{diffeq} in a neighbourhood of $s$ satisfying $f^{-1}(s) = P$, the resulting analytic map induces the map $\widehat{\kappa}$ on formal power series rings. It follows that $\widehat{\eta}$ is induced by a local period map, which completes the proof. \qed

\paragraph{Proof of \autoref{LVprop}:} ~ \\

\vspace{-0.5em}

We denote by $\mathcal{O}_{K,(v)}$ the ring of integers localized at the prime ideal $\mathfrak{p}$ of $\mathcal{O}_{K}$ corresponding to $v$. We begin by showing that (base changes to $K_{v}$ of) the points of $S(\mathcal{O}_{K,(v)})$ lie inside finitely many distinguished open affinoids $B_{p} \subset \an{S_{K_{v}}}$ admitting local period maps $\psi_{p} : B_{p} \to \an{\ch{L}_{K_{v}}}$. This reduces to showing that there are finitely many distinguished open affinoids $B_{p} \subset \an{S_{K_{v}}}$ containing the points in $S(\mathcal{O}_{K,(v)})$ over which $\an{\mathcal{H}_{K_{v}}}$ admits a rigid-analytic flat frame. We may cover $S$ by finitely many open subschemes $U \subset S$ such that $\Omega^{1}_{U}$ and the bundles $F^{k} \mathcal{H}$ for varying $k$ are all trivial. Then any point $s \in S(\mathcal{O}_{K,(v)})$ factors through some element of this cover, so we may reduce to the case where $\Omega^{1}_{S}$ and the bundles $F^{k} \mathcal{H}$ are all trivial.

Proceeding as in the proof of \autoref{padicaxschanlem}, we can choose algebraic functions $z_{1}, \hdots, z_{n}$ on $S$ such that $d z_{1}, \hdots, d z_{n}$ trivializes $\Omega^{1}_{S}$. We obtain differential equations as in (\ref{diffeqpolys}), where the polynomials $\xi_{\ell_{1}, \hdots, \ell_{e}; jk}$ are functions in the coordinate ring $R$ of $S \times \GL_{m}$, and so in particular we may view them after base-changing as elements of $R_{\mathcal{O}_{K,v}}$, where $\mathcal{O}_{K,v}$ is the ring of $v$-adic integers. Choose a point $\overline{s_{0}} \in S(\mathcal{O}_{K,N}/\mathfrak{p} \mathcal{O}_{K,N})$. Then as $S$ has (by assumption) good reduction modulo $\mathfrak{p}$, we obtain by \cite[IV. 18.5.17]{EGA} a lift $s_{0} \in S(\mathcal{O}_{K,v})$. Choosing an initial condition $P \in \GL_{m}(\mathcal{O}_{K,v})$, we may use (\ref{diffeqpolys}) to construct a map $\psi_{s_{0}} = \an{q_{K_{v}}} \circ f^{-1}$, where the partial derivatives of $f$ at $s_{0}$ are given by evaluating the polynomials $\xi_{\ell_{1}, \hdots, \ell_{e}; jk}$ at $(s_{0}, P)$. As the coefficients of the power series defining $\psi_{s_{0}}$ lie in $\mathcal{O}_{K,v}$, the map $\psi_{s_{0}}$ is defined on a residue disk $B_{p, \overline{s_{0}}}$ of radius $|p|^{1/[K_{v} : \mathbb{Q}_{p}]}$, where $|\cdot|$ is the absolute value on $\mathbb{Q}_{p}$. Varying $s_{0}$ over the finitely many elements of $S(\mathcal{O}_{K,N}/\mathfrak{p} \mathcal{O}_{K,N})$, we obtain the desired cover.

Now we wish to show that $\dim S(\mathcal{O}_{K,N})^{\textrm{Zar}} \leq d$. Recall that for each $s \in S(\mathcal{O}_{K,N})$, we have a set $\mathcal{I}(s) \subset S(\mathcal{O}_{K,N})$ of points whose associated Galois representations have isomorphic semisimplifications. As there are finitely many possibilities for the semisimplification, it suffices to consider the sets $\mathcal{S}(s) \subset S(\mathcal{O}_{K,N})$ defined by
\[ \mathcal{S}(s) = \{ s' \in S(\mathcal{O}_{K,N}) : \mathcal{I}(s) = \mathcal{I}(s') \textrm{ and } s \equiv s' \textrm{ mod } \mathfrak{p} \} , \]
and show that $\dim \overline{\mathcal{S}(s)}^{\textrm{Zar}} \leq d$. In particular, we can consider the Zariski closure of just those elements of $\mathcal{S}(s)$ whose associated points in $S(\mathcal{O}_{K,v})$ lie inside one of the neighbourhoods $B_{p,\overline{s_{0}}}$ constructed above on which we have a local period map $\psi_{s_{0}} : B_{p, \overline{s_{0}}} \to \an{\ch{L}_{K_{v}}}$. The hypothesis of the proposition tells us that the image under $\psi_{s_{0}}$ of the points in $\mathcal{S}(s)$ lie in a subvariety $O_{s} \subset \ch{L}_{\mathbb{C}_{p}}$ satisfying $\dim \Delta_{d} \geq \dim O_{s}$; here we use the fact that the flat frame on $B_{p, \overline{s_{0}}}$ is compatible with the Frobenius endomorphism (c.f. the discussion in \cite[\S3]{LV}). Base-changing to $\mathbb{C}_{p}$ and applying \autoref{padicaxschanlem} above, we find that $\dim \overline{\psi_{s_{0}}^{-1}(O_{s})}^{\textrm{Zar}} \leq d$, hence the result.

\bibliography{hodge_theory}
\bibliographystyle{alpha}

\end{document}